# Partition of a Set which Contains an Infinite Arithmetic (Respectively Geometric) Progression

Florentin Smarandache, Ph D
Professor
Department of Math & Sciences
University of New Mexico
200 College Road
Gallup, NM 87301, USA
E-mail:smarand@unm.edu

**Abstract.**
We prove that for any partition of a set which contains an infinite arithmetic (respectively geometric) progression into two subsets, at least one of these subsets contains an infinite number of triplets such that each triplet is an arithmetic (respectively geometric) progression.

**Introduction.**
First, in this article we build sets which have the following property: for any partition in two subsets, at least one of these subsets contains at least three elements in arithmetic (or geometric) progression.

**Lemma 1.** The set of natural numbers cannot be partitioned in two subsets not containing either one or the other three numbers in an arithmetic progression.
*Proof.*
Let us suppose the opposite, and have $M_1$ and $M_2$ the two subsets. Let $k \in M_1$:
a) If $k+1 \in M_1$, then $k-1$ and $k+2$ belong to $M_2$, if not we can build an arithmetic progression in $M_1$. For the same reason, since $k-1$ and $k+2$ belong to $M_2$, then $k-4$ and $k+5$ are in $M_1$. Thus $k+1$ and $k+5$ are in $M_1$ thus $k+3$ is in $M_2$; $k-4$ and $k$ are in $M_1$ thus $k+4$ is in $M_1$; we have obtained that $M_2$ contains $k+2$, $k+3$ and $k+4$, which is in contradiction with the hypothesis.
b) If $k+1 \in M_2$ then we analyze the element $k-1$. If $k-1 \in M_1$, we are in the case a) where two consecutive elements belong to the same set. If $k-1 \in M_2$, then, because $k-1$ and $k+1$ belong to $M_2$, it results that $k-3$ and $k+3 \in M_2$, then $\in M_1$. But we obtained the arithmetic progression $k-3$, $k$, $k+3$ in $M_1$, contradiction.

**Lemma 2.** If one puts aside a finite number of terms of the natural integer set, the set obtained still satisfies the property of the lemma 1.

In the lemma 1, the choice of $k$ was arbitrary, and for each $k$ one obtains at least in one of the sets $M_1$ or $M_2$ a triplet of elements in arithmetic progression: thus at least one of these two sets contains an infinity of such triplets.

If one takes a finite number of natural numbers, it takes also a finite number of triplets in arithmetic progression. But at least one of the sets $M_1$ or $M_2$ will contain an infinite number of triplets in arithmetic progression.

**Lemma 3.** If $i_1,...,i_s$ are natural numbers in arithmetic progression, and $a_1, a_2,...$ is an arithmetic (respectively geometric) progression, then $a_{i_1},....,a_{i_s}$ is also an arithmetic (respectively geometric) progression.

*Proof:*
For every $j$ we have: $2i_j = i_{j-1} + i_{j+1}$

a) If $a_1, a_2,...$ is an arithmetic progression of ratio $r$:
$$2a_{i_j} = 2(a_1 + (i_j - 1)r) = (a_1 + (i_{j-1} - 1)r) + (a_1 + (i_{j+1} - 1)r) = a_{i_{j-1}} + a_{i_{j+1}}$$

b) If $a_1, a_2,...$ is a geometric progression of ratio $r$:
$$\left(a_{i_j}\right)^2 = \left(a \cdot r^{i_j - 1}\right)^2 = a^2 \cdot r^{2i_j - 2} = \left(a \cdot r^{i_{j-1} - 1}\right) \cdot \left(a \cdot r^{i_{j+1} - 1}\right) = a_{i_{j-1}} + a_{i_{j+1}}$$

**Theorem 1.**
It does not matter the way in which one partitions the set of the terms of an infinite arithmetic (respectively geometric) progression in subsets: in at least one of these subsets there will be at least three terms in arithmetic (respectively geometric) progression.

*Proof:*
According to lemma 3, it is enough to study the partition of the set of the indices of the terms of the progression in 2 subsets, and to analyze the existence (or not) of at least 3 indices in arithmetic progression in one of these subsets.
But the set of the indices of the terms of the progression is the set of the natural numbers, and we proved in lemma 1 that it cannot be partition in 2 subsets without having at least 3 numbers in arithmetic progression in one of these subsets: the theorem is proved.

**Theorem 2.**
A set $M$, which contains an infinite arithmetic (respectively geometric) not constant progression, preserves the property of the theorem 1.

Indeed, this directly results from the fact that any partition of $M$ implies the partition of the terms of the progression.

**Application:** Whatever is the way in which one partitions the set $A = \{1^m, 2^m, 3^m,...\}$, $(m \in N^*)$ in subsets, at least one of these subsets contains three terms in geometric progression.

(Generalization of the problem 0:255 from "Gazeta Matematică", Bucharest, No. 10/1981, p. 400).

The solution naturally results from theorem 2, if it is noticed that $A$ contains the geometric progression $a_n = (2^m)^n$, $(n \in \mathbb{N}^*)$.

Moreover one can prove that in at least one of the subsets there is an infinity of triplets in geometric progression, because $A$ contains an infinity of different geometric progressions: $a_n^{(p)} = (p^m)^n$ with $p$ prime and $n \in \mathbb{N}^*$, to which one can apply the theorems 1 and 2.

**Reference**:

F. Smarandache, "About Some Progressions", in Collected Papers, Vol. I, Ed. Tempus, Bucharest, pp. 60-62, 1996.